\documentclass[reqno]{amsart}
\usepackage[active]{srcltx}
\usepackage{amsmath}
\usepackage{times}
\usepackage[T1]{fontenc}
\usepackage{textcomp}
\usepackage{txfonts}
\usepackage{tikz}
\usetikzlibrary{arrows,decorations.markings}
\usetikzlibrary{shapes}
\usepackage{mathrsfs}
\usepackage{multicol}
\usepackage{hyperref}
\usepackage{amsthm}

\hypersetup{pdftitle={  TITLE },
  pdfauthor={Ge Li},
  colorlinks = true,
  linkcolor =blue,
  anchorcolor = red,
  citecolor = blue,
  urlcolor = blue
}

\usepackage{titlesec}

\titleformat{\section}{\Large\bfseries}{\thesection.}{0.5em}{}
\titleformat{\subsection}{\large\bfseries}{\thesubsection}{0.5em}{}

\usepackage{pifont}

\usepackage{bm}

\DeclareMathAlphabet{\mathpzc}{OT1}{pzc}{m}{it}

{\catcode`\|=\active
  \gdef\set#1{\mathinner{\lbrace\,{\mathcode`\|"8000%
                                   \let|\midvert #1}\,\rbrace}}
}
\def\midvert{\egroup\mid\bgroup}

\usepackage{aliascnt}
\def\NewTheorem#1{%
  \newaliascnt{#1}{equation}
  \newtheorem{#1}[#1]{#1}
  \aliascntresetthe{#1}
  \expandafter\def\csname #1autorefname\endcsname{#1}
}

\newcounter{main}
\theoremstyle{plain}

\numberwithin{equation}{section}
\NewTheorem{Assumption}
\NewTheorem{Proposition}
\NewTheorem{Theorem}
\NewTheorem{Corollary}
\NewTheorem{Conjecture}
\NewTheorem{Lemma}
\NewTheorem{Definition}
\theoremstyle{definition}
\NewTheorem{Remark}
\NewTheorem{Remarks}
\NewTheorem{Example}

\def\@thm#1#2#3{%
  \ifhmode\unskip\unskip\par\fi
  \normalfont
  \trivlist
  \let\thmheadnl\relax
  \let\thm@swap\@gobble
  \let\thm@indent\noindent
  \thm@headfont{\SMStheoremfont}
  \thm@notefont{\fontseries\mddefault\upshape}%
  \thm@headpunct{.}
  \thm@headsep 5\p@ plus\p@ minus\p@\relax
  \thm@space@setup
  #1
  \@topsep \thm@preskip               
  \@topsepadd \thm@postskip           
  \def\@tempa{#2}\ifx\@empty\@tempa
    \def\@tempa{\@oparg{\@begintheorem{#3}{}}[]}%
  \else
    \refstepcounter{#2}%
    \def\@tempa{\@oparg{\@begintheorem{#3}{\csname the#2\endcsname}}[]}%
  \fi
  \@tempa
}

\newcommand{\Z}{\mathbb{Z}}
\newcommand{\N}{\mathbb{N}}

\newcommand{\F}{\mathbb{F}}
\newcommand{\Lr}{L_r(\bi)^{\Nr}}
\newcommand{\Nr}{N_r(\bi)}

\def\bi{\mathbf{i}}
\def\bj{\mathbf{j}}
\def\t{\mathsf{t}}

\def\l{\ell}

\usepackage{tikz}  

\colorlet{darkgreen}{green!50!black}
\tikzset{dots/.style={ultra thick,loosely dotted},
         belt/.style={draw,blue,thick,fill=blue!50},
         greendot/.style={fill,circle,color=darkgreen,inner sep=1.5pt,outer sep=0}
}
\newenvironment{braid}{
  \begin{tikzpicture}[baseline=6mm,blue,line width=1pt, xscale = 0.35, yscale=0.4,
                      draw/.append style={rounded corners},
                      every node/.append style={font=\fontsize{5}{5}\selectfont}]%
  }{\end{tikzpicture}
}

\def\Grid(#1,#2){
  \draw[very thin,gray,step=2mm] (0,0)grid(#1,#2);
  \draw[very thin,darkgreen,step=10mm] (0,0)grid(#1,#2);
}

\newcommand\Tableau[2][\relax]{
  \begin{tikzpicture}[scale=0.5,draw/.append style={thick,black},baseline=-3mm]
    \ifx\relax#1\relax%
    \else 
      \foreach\box in {#1} {
        \ifx\box\relax\else
          \filldraw[blue!10]\box+(-.5,-.5) rectangle ++(.5,.5);
        \fi
      }
    \fi;
    \newcount\row\newcount\col
    \row=0
    \foreach \Row in {#2} {
       \col=1
       \foreach \k in \Row {
          \draw (\the\col,\the\row) +(-.5,-.5) rectangle ++(.5,.5);
          \ifnum\k<0
            \draw(\the\col,\the\row) node[fill=gray!20]{-\k};
          \else \draw (\the\col,\the\row) node{\k};
          \fi
          \global\advance\col by 1
       }
       \global\advance\row by -1
    }
  \end{tikzpicture}
}
\def\Tritab(#1|#2|#3){\Bigg(\hspace*{1mm}\Tableau{#1}\hspace*{1mm}\Bigg|\hspace*{1mm}%
    \Tableau{#2}\hspace*{1mm}\Bigg|\hspace*{1mm}\Tableau{#3}\hspace*{1mm}\Bigg)
}

\newcommand\YoungDiagram[2][\relax]{
  \begin{tikzpicture}[scale=0.5,draw/.append style={thick,black},baseline=-1mm]
    \ifx\relax#1\relax%
    \else 
    \foreach\box in {#1} {
      \filldraw[blue!10]\box rectangle ++(1,1);
    }
    \fi
    \newcount\row
    \row=0
    \foreach \col in {#2} {
       \draw(1,\the\row)grid ++(\col,1);
       \global\advance\row by -1
    }
  \end{tikzpicture}
}

\def\Tridiag(#1|#2|#3){\Bigg(\hspace*{1mm}\YoungDiagram{#1}\hspace*{1mm}\Bigg|\hspace*{1mm}%
    \YoungDiagram{#2}\hspace*{1mm}\Bigg|\hspace*{1mm}\YoungDiagram{#3}\hspace*{1mm}\Bigg)
}

\DeclareMathOperator{\res}{res}

\newcommand{\map}[2]{\,{:}\,#1\!\longrightarrow\!#2}

\usepackage[enableskew,vcentermath]{youngtab}
\def\tab(#1){\mbox{\tiny$\young(#1)$}\,}
\def\ydiag(#1){\mbox{\tiny$\yng(#1)$}\,}

\newpage

\begin{document}    
\bibliographystyle{andrew}

\title[Periodic Property]{\larger[3] The idempotents in cyclotomic Hecke algebras and periodic property of the Jucys-Murphy elements\\[40pt]}
\author{\larger[2]Ge~Li}
\address{School of Mathematics and Statistics\\
		University of Sydney\\
		Sydney, NSW 2006}
\email{geli@maths.usyd.edu.au}
\begin{abstract}
This paper proves a periodic property of Jucys-Murphy elements of the degenerate and non-degenerate cyclotomic Hecke algebras of type A. We do this by first giving a new closed formula for the KLR idempotents $e(\bi)$ which, it tuns out, is very efficient computationally.
\end{abstract}

\maketitle          
\pagenumbering{arabic}

\section{Introduction}

The degenerate and non-degenerate cyclotomic Hecke algebras $H_n^\Lambda$ of type A are important algebras because they arise in the categorification of the canonical basis of the affine special linear groups~\cite{BK:GradedDecomp}. Recently, building on works of Khovanov and Lauda~\cite{KhovLaud:diagII, KhovLaud:diagI} and Rouquier~\cite{Rouq:2KM}, Brundan and Kleshchev~\cite{BK:GradedKL} showed that $H_n^\Lambda$ is isomorphic to a Khovanov-Lauda-Rouquier algebra of type A, and hence that these algebras are $\Z$-graded.

Central to construction of Brundan and Kleshchev graded isomorphism theorem are certain deponents $e(\bi)\in H_n^\Lambda$, which in the case of the symmetric groups have their origins in the work of Murphy~\cite[(1,2)]{M:Nak}. The purpose of this paper is to give explicit formulas for these idempotents when $e > 0$ and $p > 0$, which are very efficient computationally, and to use them to prove a periodicity property for the Jucys-Murphy elements of $H_n^\Lambda$.

This paper is organized as follows. In Section 2 we introduced the definition of (degenerate and non-degenerate) cyclotomic Hecke algebras of type A. We then give the definition of $e(\bi)$'s constructed by Kleshchev~\cite[Lemma 7.1]{Klesh:book} in $H_n^\Lambda$ and a set of nilpotency elements $\set{y_r|1\leq r\leq n}$ defined by Brundan and Kleshchev~\cite[(3.21),(4.21)]{BK:GradedKL}. In Section 3 we gave an explicit expression of idempotents $e(\bi)$. In Section 4 and 5 we proved the periodic property of Jucys-Murphy elements in degenerate and non-degenerate cyclotomic Hecke algebras of type A, respectively. In Section 6 we improved our results of Section 3 and showed that an expression of $e(\bi)$'s depended on the nilpotency degree of $y_r$'s.

This paper is a part of my PhD thesis~\cite{Li:thesis}.

\textit{Acknowledgements.} We thank Andrew Mathas for his support while writing this paper.

\section{Cyclotomic Hecke algebras of type A and the idempotents $e(\bi)$}

Let $\F_p$ be a fixed field of characteristic $p \geq 0$ with $q\in \F_p^\times$. Let $e$ be the smallest positive integer such that $1 + q + \ldots + q^{e-1} = 0$ and setting $e = 0$ if no such integer exists. Then define $I = \Z/e\Z$ if $e > 0$ and $I = \Z$ if $e = 0$. 

For $n \geq 0$, assume that $q = 1$. Let $H_n$ be the \textbf{degenerate affine Hecke algebra}, working over $\F_p$. So $H_n$ has generators
$$
\{x_1,\ldots,x_n\}\cup\{s_1,\ldots,s_{n-1}\}
$$
subject to the following relations
$$
\begin{array}{cc}
x_r x_s = x_s x_r; & \\
s_r x_{r+1} = x_r s_r + 1, & s_r s_x = x_s s_r \hspace*{5mm}\text{if $s \neq r,r+1$;}\\
s_r^2 = 1; \\
s_r s_{r+1}s_r = s_{r+1}s_r s_{r+1}, & s_r s_t = s_t s_r \hspace*{5mm}\text{if $|r - t| > 1$.}
\end{array}
$$

Now we assume that $q\neq 1$ and $H_n$ be the \textbf{non-degenerate affine Hecke algebra} over $\F_p$. So $H_n$ has generators
$$
\{ X_1^{\pm 1},\ldots,X_n^{\pm 1}\}\cup\{T_1,\ldots,T_{n-1}\}
$$
subject to the following relations
$$
\begin{array}{cc}
X_r^{\pm 1} X_s^{\pm 1} = X_s^{\pm 1} X_r^{\pm1}, & X_r X_r^{-1} = 1;\\
T_r X_r T_r = q X_{r+1}, & T_r X_s = X_s T_r\hspace*{5mm}\text{if $s\neq r,r+1$};\\
T_r^2 = (q-1)T_r + q; &\\
T_r T_{r+1}T_r = T_{r+1}T_rT_{r+1},& T_r T_s = T_s T_r\hspace*{5mm}\text{if $|r-s| > 1$}.
\end{array}
$$

To the index set $I$, we associate two lattices
$$
P:= \bigoplus_{i\in I}\Z \Lambda_i, \qquad Q:= \bigoplus_{i\in I}\Z\alpha_i,
$$
and let $(.,. )\map{P\times Q}{\Z}$ be the bilinear pairing defined by $(\Lambda_i,\alpha_j) = \delta_{i,j}$. Let $P_+$ and $Q_+$ denote the subset of $P$ and $Q$ such that $P_+ = \bigoplus_{i\in I}\N \Lambda_i$ and $Q_+ = \bigoplus_{i\in I}\N \alpha_i$.

Then for any $\Lambda \in P_+$, we define
\begin{equation} \label{E:QHecke}
H_n^\Lambda = \begin{cases}
H_n/ \langle \prod_{i\in I}(X_1 - q^i)^{(\Lambda,\alpha_i)}\rangle, & \text{if $q\neq 1$;}\\
H_n/ \langle \prod_{i\in I}(x_1 - i)^{(\Lambda,\alpha_i)}\rangle, & \text{if $q = 1$.}
\end{cases}
\end{equation}
and we call $H_n^\Lambda$ the \textbf{degenerate cyclotomic Hecke algebra} if $q = 1$ and \textbf{non-degenerate cyclotomic Hecke algebra} if $q\neq 1$. We call $p$ the \textbf{characteristic} of $H_n^\Lambda$ and $e$ the \textbf{quantum characteristic} of $H_n^\Lambda$. 

By the definitions, degenerate and non-degenerate cyclotomic Hecke algebras are similar with some minor difference. In order to minimize their difference we define 
\begin{equation} \label{E:q_i}
q_i = \begin{cases}
i, & \text{if $q = 1$;}\\
q^i, & \text{if $q\neq 1$.}
\end{cases}
\end{equation}
and use $x_r$ instead of $X_r$ when we don't have to distinguish which case we are working with. Hence we can rewrite (\ref{E:QHecke}) as
\begin{equation} \label{E:cycloHecke}
H_n^\Lambda = H_n/\langle \prod_{i\in I}(x_1 - q_i)^{(\Lambda,\alpha_i)}\rangle.
\end{equation}

In $H_n^\Lambda$, the elements $x_r$ and $X_r$ are called the \textbf{Jucys-Murphy elements}. 

We can define a set of pairwise orthogonal idempotents $\set{e(\bi)|\bi\in I^n}$ for both degenerate and non-degenerate $H_n^\Lambda$. By Kleshchev~\cite[Lemma 7.1]{Klesh:book}, the eigenvalues of each $x_r$ or $X_r$ on $M$ belongs to $I$. So $M = \bigoplus_{\bi\in I^n} M_\bi$ of its weight space
\begin{equation} \label{E:e(i)}
M_\bi = \set{v\in M | (x_r - q_{i_r})^N v = 0\text{ for all $r = 1,\ldots,d$ and $N\gg 0$}} ,
\end{equation}
where $q_{i_r}$ is introduced in (\ref{E:q_i}). Then we deduce that there is a system $\set{e(\bi)|\bi\in I^n}$ such that $e(\bi)M = M_\bi$, and $\bi \in I^n$ is called the residue sequence.

Brundan and Kleshchev \cite[(3.21),(4.21)]{BK:GradedKL} introduced a set of elements $y_r$ in both degenerate and non-degenerate cyclotomic Hecke algebras, where
\begin{equation}
y_r := \begin{cases} \label{E:y_r}
\sum_{\bi\in I^n} (x_r - i_r)e(\bi), & \text{if $q = 1$;}\\
\sum_{\bi\in I^n} (1 - q^{-i_r}X_r) e(\bi). & \text{if $q \neq 1$.}
\end{cases}
\end{equation}
and by~\cite[Lemma 2.1]{BK:GradedKL}, in both cases $y_r^s = 0$ for $s\gg 0$. It is easy to imply that $(X_r - q^{i_r})^s e(\bi) = 0$ for $s \gg 0$ in non-degenerate case. Hence we have the following Corollary.

\begin{Corollary} \label{nil}
We have $(x_r - q_{i_r})^s e(\bi) = 0$ for $s \gg 0$ in both degenerate and non-degenerate cyclotomic Hecke algebras. Furthermore, $(x_r - q_{i_r})^s e(\bi) = 0$ if and only if $y_r^s e(\bi) = 0$. 
\end{Corollary}

We will apply \autoref{nil} without mention.

The smallest $s$ with $y_r^s e(\bi) = 0$ is called the \textbf{nilpotency degree} of $y_r$ corresponds to $e(\bi)$. Using the explicit expression of $e(\bi)$ we found, we can find the periodic property of the Jucys-Murphy elements in both degenerate and non-degenerate cyclotomic Hecke algebras when $e > 0$ and $p > 0$, and the period of $x_r$ is determined by the nilpotency degree of $y_r$'s.

\section{Explicit expression of idempotents $e(\bi)$ part I}

\subsection{General cases}

In this section, first we introduce some properties of characteristic and quantum characteristic.

\begin{Lemma} ~\label{useful}
Suppose $\F_p$ is a field with $\text{char }\F_p = p > 0$ and $r_1, r_2 \in H_n^\Lambda$. For any non-negative integer $k$ we always have $(r_1 - r_2)^{p^k} = r_1^{p^k} - r_2^{p^k}$.
\end{Lemma}

\begin{Remark}
Notice the above Lemma is a well-known result and will be applied without mention in this chapter.
\end{Remark}

\begin{Lemma} \label{e,p in deg}
Suppose $p$ and $e$ are characteristic and quantum characteristic of degenerate $H_n^\Lambda$. Then $e = p$.
\end{Lemma}

\proof The Lemma follows directly by the definitions of $p$ and $e$ in degenerate $H_n^\Lambda$. \endproof

\begin{Lemma} \label{e,p in nondeg}
Suppose $p$ and $e$ are characteristic and quantum characteristic of non-degenerate $H_n^\Lambda$ with $e,p > 0$. Then $\text{gcd}(e,p) = 1$. Moreover, we can find $l$ such that $p^l \equiv 1 \pmod{e}$.
\end{Lemma}

\proof In non-degenerate case, $\text{gcd}(e,p) = 1$ is well-known. So by Chinese Remainder Theorem we can find $a,b\in \Z$ such that $ap + be = 1$. Now consider the sequence $p,p^2,p^3,p^4\ldots$. We can find $k_1,k_2$ such that $p^{k_1} \equiv p^{k_2} \pmod{e}$. Choose $k_2$ such that $k_2 - k_1 > k_1$. Hence write $l = k_2 - k_1$ and $p^l \equiv s \pmod{e}$ where $0 \leq s \leq e-1$. So we have $p^{2l} \equiv p^l \pmod{e}$ which implies $s^2 \equiv s \pmod{e}$. So we can write $s^2 - s = ke$ for some $k \in \Z$. So
$$
s^2 - s = ke \hspace*{3mm} \Rightarrow \hspace*{3mm} as(s-1) = ake \hspace*{3mm} \Rightarrow \hspace*{3mm} (1 - be)(s-1) = ake \hspace*{3mm} \Rightarrow \hspace*{3mm} s - 1 = (b(s-1) + ak)e
$$
which implies $e\ |\ s-1$. But because $0 \leq s \leq e-1$, we have $s = 1$. Therefore $p^l \equiv 1\pmod{e}$. \endproof

Therefore, in degenerate case, by \autoref{e,p in deg} we have $e = p$ and in non-degenerate case, by \autoref{e,p in nondeg} we have $\text{gcd}(e,p) = 1$, i.e. $e \neq p$.

In the rest of this paper we set $p > 0$ and $e> 0$.  Fix a residue sequence $\bi = (i_1,i_2,\ldots,i_n) \in I^n$. For any $1\leq r\leq n$ and any $j\in I$ with $j\neq i_r$, choose $N \gg 0$ and define $L_{i_r,j} = 1 - (\frac{q_{i_r} - x_r}{q_{i_r} - q_j})^N$ in both degenerate and non-degenerate cases.

Notice that by the definition of $e(\bi)$ given in (\ref{E:e(i)}), for any $\bj\in I$ and $1 \leq r \leq n$, we have
$$
(x_r - q_{j_r})^N e(\bj) = 0
$$
for $N \gg 0$.

\begin{Lemma} \label{deg:1}
Suppose $1\leq r\leq n$ and $\bj = (j_1,j_2,\ldots,j_n) \in I^n$. For $j\in I$ and $N_j \gg 0$ we have
$$
L_{i_r,j}^{N_j}e(\bj) = \begin{cases}
e(\bj), & \text{if $j_r = i_r$};\\
0, & \text{if $j_r = j$}.
\end{cases}
$$
\end{Lemma}

\proof Suppose $j_r = i_r$. Because $(x_r - q_{i_r})^N e(\bj) = (x_r - q_{j_r})^N e(\bj) = 0$ for $N \gg 0$, we have
$$
L_{i_r,j}e(\bj) = (1 - (\frac{q_{i_r} - x_r}{q_{i_r} - q_j})^N)e(\bj) = e(\bj) - \frac{1}{(q_{i_r} - q_j)^N}(q_{i_r} - x_r)^Ne(\bj) = e(\bj).
$$

Therefore $L_{i_r,j}^{N_j}e(\bj) = L_{i_r,j}^{N_j - 1} e(\bj) = \ldots = L_{i_r,j}e(\bj) = e(\bj)$.

Suppose $j_r = j$. We have
\begin{eqnarray*}
L_{i_r,j} & = &1 - (\frac{q_{i_r} - x_r}{q_{i_r} - q_j})^N\\
& = & -\sum_{k = 1}^{N} (\frac{q_{i_r} - x_r}{q_{i_r} - q_j})^k + \sum_{k = 0}^{N-1} (\frac{q_{i_r} - x_r}{q_{i_r} - q_j})^k\\
& = & \frac{x_r - q_{i_r}}{q_{i_r} - q_j}\sum_{k = 0}^{N-1} (\frac{q_{i_r} - x_r}{q_{i_r} - q_j})^k + \frac{q_{i_r} - q_j}{q_{i_r} - q_j}\sum_{k = 0}^{N-1} (\frac{q_{i_r} - x_r}{q_{i_r} - q_j})^k\\
& = & (x_r - q_{i_r} + q_{i_r} - q_j)\frac{1}{q_{i_r} - q_j}\sum_{k = 0}^{N-1} (\frac{q_{i_r} - x_r}{q_{i_r} - q_j})^k\\
& = & (\frac{1}{q_{i_r} - q_j}\sum_{k = 0}^{N-1} (\frac{q_{i_r} - x_r}{q_{i_r} - q_j})^k) (x_r - q_j).
\end{eqnarray*}

Therefore for $N_j \gg 0$,
$$
L_{i_r,j}^{N_j}e(\bj) = (\frac{1}{q_{i_r} - q_j}\sum_{k = 0}^{N-1} (\frac{q_{i_r} - x_r}{q_{i_r} - q_j})^k)^{N_j} (x_r - q_j)^{N_j} e(\bj) = 0
$$
because when $j_r = j$ we have $(x_r - q_j)^{N_j} e(\bj) = 0$, which completes the proof. \endproof

Now we define $L_r(\bi) = \prod_{\substack{j\in I\\j\neq i_r}} L_{i_r,j}$. In the product, $j \in I\backslash\{i_r\}$, which is a finite product since $e > 0$. So $L_r(\bi)$ is well defined. We have the following Lemma.

\begin{Lemma} \label{deg:2}
Suppose $1\leq r\leq n$. We can choose $\Nr \gg 0$ such that
$$
\Lr  = \sum_{\substack{\bj \in I^n \\ j_r = i_r}} e(\bj).
$$
\end{Lemma}

\proof By \autoref{deg:1}, for any $j\in I$ with $j \neq i_r$ we can find $N_j$ large enough such that
$$
L_{i_r,j}^{N_j}e(\bj) = \begin{cases}
e(\bj), & \text{if $j_r = i_r$},\\
0, & \text{if $j_r = j$}.
\end{cases}
$$
Now choose $\Nr \geq \max\{N_j\ |\ j\in I, j \neq i_r\}$, which is finite since $e > 0$. Therefore $\Lr  = \prod_{\substack{j\in I\\j\neq i_r}} L_{i_r,j}^{\Nr}$. Hence for any $e(\bj)$, if $j_r \neq i_r$, 
\begin{equation} \label{deg:2:E1}
\Lr e(\bj) = \prod_{\substack{j\in I\\j\neq i_r}} L_{i_r,j}^{\Nr}e(\bj) = (\prod_{\substack{j\in I\\j\neq i_r,j_r}} L_{i_r,j}^{\Nr})L_{r,j_r}^{\Nr} e(\bj) = 0,
\end{equation}
and if $j_r = i_r$,
\begin{equation} \label{deg:2:E2}
\Lr e(\bj) = e(\bj),
\end{equation}
because for any $j$, $L_{i_r,j}^{N_j}e(\bj) = e(\bj)$.

Therefore, because $\sum_{\bj\in I^n} e(\bj) = 1$, by (\ref{deg:2:E1}) and (\ref{deg:2:E2}), 
$$
\Lr  = \Lr (\sum_{\bj\in I^n} e(\bj)) = \sum_{\bj \in I^n}\Lr e(\bj) = \sum_{\substack{\bj \in I^n \\ j_r = i_r}}\Lr  e(\bj) = \sum_{\substack{\bj \in I^n \\ j_r = i_r}} e(\bj)
$$
which completes the proof. \endproof

As the idempotents $e(\bj)$'s are pairwise orthogonal, \autoref{deg:2} immediately implies the following.

\begin{Corollary} \label{deg:3}
For any $\bi \in I^n$, we have
$$
e(\bi) = \prod_{r = 1}^n \Lr.
$$
\end{Corollary}

The previous results are true in both degenerate and non-degenerate cases. Notice that when we define $L_{i_r,j} = 1 - (\frac{q_{i_r} - x_r}{q_{i_r} - q_j})^{N_j}$ and $\Lr $, the only restriction is that $N_j$ and $\Nr$ are large enough. As we now show, by choosing specific values for $N_j$ and $\Nr$, it is possible to simplify the expression of $ \Lr $ even further and give a more explicit expression of $e(\bi)$. We emphasize the simplified expressions for $ \Lr $ are different for degenerate and non-degenerate $H_n^\Lambda$.

\subsection{Explicit expression in degenerate cases}
We start with the degenerate cyclotomic Hecke algebras. Recall that in this case $e = p$.

\begin{Proposition} \label{deg:4}
Suppose $q = 1$. For any $i_r \in I$ there exists $s \gg 0$ such that
$$
\sum_{\substack{\bj \in I^n \\ j_r = i_r}}e(\bj) = \begin{cases}
1 - x_r^{p^s(1-p)}, & \text{when $i_r = 0$},\\
-\sum_{k = 1}^{p-1} \frac{x_r^{kp^s}}{i_r^k}, & \text{when $i_r \neq 0$}.
\end{cases}
$$
\end{Proposition}

\proof By \autoref{deg:2} the Proposition is equivalent to claim that
$$
\Lr  = \begin{cases}
1 - x_r^{p^s(1-p)}, & \text{when $i_r = 0$},\\
-\sum_{k = 1}^{p-1} \frac{x_r^{kp^s}}{i_r^k}, & \text{when $i_r \neq 0$}.
\end{cases}
$$
for $s \gg 0$.

By the definition of $L_r(\bi)$, because $I = \Z/p\Z$ we have
$$
L_r(\bi) = \prod_{\substack{j\in I\\j\neq i_r}} L_{i_r,j} = \prod_{\substack{j\in I\\j\neq i_r}} (1 - (\frac{i_r - x_r}{i_r - j})^{N_j}) = \prod_{j = 1}^{p-1}(1 - (\frac{i_r - x_r}{j})^{N_j}).
$$

Take $k \gg 0$ and $N_j = p^k$. Hence because $H_n^\Lambda$ is defined over a field $\F_p$ of characteristic $p$, we have $j^{N_j} = j$. And because $p$ is a prime, we have
$$
L_r(\bi) = \prod_{j = 1}^{p-1}(1 - (\frac{i_r - x_r}{j})^{N_j}) = \prod_{j = 1}^{p-1}(1 - \frac{(i_r - x_r)^{N_j}}{j}) = \prod_{j = 1}^{p-1}(1 - j{\cdot}(i_r - x_r)^{N_j}) = 1 - (i_r - x_r)^{(p-1)N_j}.
$$

Without loss of generality, choose $\Nr = p^l$ with $l\gg 0$. We have
$$
\Lr  = (1 - (i_r - x_r)^{(p-1)N_j})^{\Nr} = 1 - (i_r - x_r)^{(p-1)N_j\Nr}.
$$

Setting $s = k+l$, we have $N_j\Nr = p^{k+l} = p^s$. Now we consider two cases, which are $i_r = 0$ and $i_r \neq 0$.

Suppose first $i_r = 0$. We have
\begin{equation} \label{deg:4:E1}
\Lr  = 1 - (i_r - x_r)^{(p-1)N_j\Nr} = 1 - (-x_r)^{(p-1)p^s} = 1 - x_r^{(p-1)p^s}.
\end{equation}

Suppose $i_r \neq 0$. We have
\begin{eqnarray*}
(i_r - x_r)^{(p-1)N_j\Nr} & = & (i_r - x_r)^{p^{s+1} - p^s} = \frac{(i_r - x_r)^{p^{s+1}}}{(i_r - x_r)^{p^s}}\\
& = & \frac{ i_r - x_r^{p^{s+1}}}{i_r - x_r^{p^s}} = \frac{1 - (\frac{x_r}{i_r})^{p^{s+1}}}{1 - (\frac{x_r}{i_r})^{p^s}}\\
& = & 1 + (\frac{x_r}{i_r})^{p^s} + (\frac{x_r}{i_r})^{2p^s} + \ldots + (\frac{x_r}{i_r})^{(p-1)p^s}\\
& = & 1 + \frac{x_r^{p^s}}{i_r} + \frac{x_r^{2p^s}}{i_r^2} + \ldots + \frac{x_r^{(p-1)p^s}}{i_r^{p-1}} = \sum_{k = 0}^{p-1} \frac{x^{kp^s}}{i_r^k}.
\end{eqnarray*}

Hence,
\begin{equation} \label{deg:4:E2}
\Lr  = 1 - (i_r - x_r)^{(p-1)N_j\Nr} = 1 - \sum_{k = 0}^{p-1} \frac{x^{kp^s}}{i_r^k} = -\sum_{k = 1}^{p-1} \frac{x^{kp^s}}{i_r^k}.
\end{equation}

By combining (\ref{deg:4:E1}) and (\ref{deg:4:E2}), we complete the proof. \endproof

Finally, by combining \autoref{deg:3} and \autoref{deg:4}, we have an explicit expression of $e(\bi)$ for the degenerate cyclotomic Hecke algebras.

\begin{Theorem} \label{deg:exp1}

Suppose $\bi = (i_1,i_2,\ldots,i_n) \in I^n$ and $q = 1$, then
$$
e(\bi) = \prod_{r = 1}^n \Lr 
$$
where
$$
\Lr  = \begin{cases}
1 - x_r^{p^s(1-p)}, & \text{when $i_r = 0$},\\
-\sum_{k = 1}^{p-1} \frac{x_r^{kp^s}}{i_r^k}, & \text{when $i_r \neq 0$}.
\end{cases}
$$
for $s \gg 0$.

\end{Theorem}

\subsection{Explicit expression in non-degenerate cases}
We now give a similar expression for the non-degenerate cyclotomic Hecke algebras. First we give two Lemmas which will be used later.

\begin{Lemma} \label{nondeg:1}
For any $k\in \Z$ with $k\not\equiv 0\pmod{e}$, we have
$$
1 + q^k + q^{2k} + \ldots + q^{(e-1)k} = 0.
$$
\end{Lemma}

\proof By the definition, we have
\begin{eqnarray*}
&& 1 + q + q^2 + \ldots + q^{e-1} = 0\\
& \Rightarrow & (1 + q + q^2 + \ldots + q^{e-1})(1-q) = 0\\
& \Rightarrow & 1 - q^e = 0\\
& \Rightarrow & q^e = 1\\
& \Rightarrow & (q^e)^k = q^{ke} = (q^k)^e = 1\\
& \Rightarrow & (q^k)^e - 1 = (1 + q^k + q^{2k} + \ldots + q^{(e-1)k})(q^k - 1) = 0.
\end{eqnarray*}

Because $k\in \Z$ and $k\not\equiv 0\pmod{e}$, we have $q^k - 1 \neq 0$. Therefore we must have $1 + q^k + q^{2k} + \ldots + q^{(e-1)k} = 0$. \endproof

\begin{Lemma} \label{nondeg:2}
Suppose $i_r \in I$ and $f(x) = \prod_{j\neq i_r}(1 - \frac{r^{i_r} - x}{r^{i_r} - r^j}) \in \F_p[x]$ with $r = q^s$ for some positive integer $s \not\equiv 0\pmod{e}$ and $q \in \F_p^{\times}$. Then $e^{-1} \in \F_p$ and
$$
f(x) = e^{-1}(1 + \frac{x}{r^{i_r}} + (\frac{x}{r^{i_r}})^2 + \ldots + (\frac{x}{r^{i_r}})^{e-1}).
$$
\end{Lemma}

\proof By \autoref{e,p in nondeg} we have $\text{gcd}(e,p) = 1$ and hence $e^{-1} \in \F_p$. Define $g(x) = e^{-1}(1 + \frac{x}{r^{i_r}} + (\frac{x}{r^{i_r}})^2 + \ldots + (\frac{x}{r^{i_r}})^{e-1})$. We prove that $f(x) = g(x)$ by first comparing their roots. It is obvious that the roots of $f(x)$ are all of the form $r^j$ with $j\in I$ and $j \neq i_r$. Then for any such $r^j$,
$$
g(r^j) = e^{-1}(1 + r^{j - i_r} + r^{2(j - i_r)} + \ldots + r^{(e-1)(j - i_r)}) =  e^{-1}(1 + r^k + r^{2k} + \ldots + r^{(e-1)k})
$$
for $k \equiv j - i_r \pmod{e}$ and $k\neq 0$. Because $r = q^s$ and $s \not\equiv 0\pmod{e}$, we must have $sk \not\equiv 0\pmod{e}$. Therefore by \autoref{nondeg:1} we have $g(r^j) = 0$. Because $f(x)$ and $g(x)$ are both polynomials of degree $e-1$, they have $e-1$ roots, which means that $g(x)$ and $f(x)$ have the same roots. This yields that $f(x) = kg(x)$ for some $k\in \F_p$.

Now because $f(r^{i_r}) = 1 = g(r^{i_r})$, we have $k = 1$. Therefore $f(x) = g(x)$, which completes the proof. \endproof

\begin{Proposition} \label{nondeg:3}
Suppose $q \neq 1$. For any $i_r \in I$, there exists $s \gg 0$ such that
$$
\sum_{\substack{\bj \in I^n \\ j_r = i_r}}e(\bj) = e^{-1}(1 + \frac{X_r^{p^s}}{q^{p^s{\cdot}i_r}} + (\frac{X_r^{p^s}}{q^{p^s{\cdot}i_r}})^2 + \ldots + (\frac{X_r^{p^s}}{q^{p^s{\cdot}i_r}})^{e-1}).
$$
\end{Proposition}

\proof By \autoref{deg:2} the Proposition is equivalent to prove that
$$
\Lr = e^{-1}(1 + \frac{X_r^{p^s}}{q^{p^s{\cdot}i_r}} + (\frac{X_r^{p^s}}{q^{p^s{\cdot}i_r}})^2 + \ldots + (\frac{X_r^{p^s}}{q^{p^s{\cdot}i_r}})^{e-1}).
$$

By the definition of $L_r(\bi)$, because $I = \Z/e\Z$, if $N_j,\Nr\gg 0$ then we have
$$
L_r(\bi)^{\Nr} = \prod_{\substack{j\in I\\ j\neq i_r}}(1 - (\frac{q^{i_r} - X_r}{q^{i_r} - q^j})^{N_j})^{\Nr}.
$$

Suppose $N_j = p^k$ and $\Nr = p^l$ with $k,l\gg 0$. We have
\begin{eqnarray*}
L_r(\bi)^{\Nr} & = & \prod_{j\neq i_r}(1 - (\frac{q^{i_r} - X_r}{q^{i_r} - q^j})^{p^k})^{p^l}\\
& = & \prod_{j\neq i_r}(1 - (\frac{q^{i_r} - X_r}{q^{i_r} - q^j})^{p^{k+l}})\\
& = & \prod_{j\neq i_r}(1 - \frac{q^{p^{k+l}{\cdot}i_r} - X_r^{p^{k+l}}}{q^{p^{k+l}{\cdot}i_r} - q^{p^{k+l}{\cdot}j}})\\
& = & \prod_{j\neq i_r}(1 - \frac{r^{i_r} - X_r^{p^{s}}}{r^{i_r} - r^j}),
\end{eqnarray*}
where $s = k+l$ and $r = q^{p^s}\in \F_p$. Notice that by \autoref{e,p in nondeg}, we have $p^s \not\equiv 0\pmod{e}$.

Now we set $f(x) =  \prod_{j\neq i_r}(1 - \frac{r^{i_r} - x}{r^{i_r} - r^j}) \in \F_p[x]$. By \autoref{nondeg:2} we have
$$
f(x) = e^{-1}(1 + \frac{x}{r^{i_r}} + (\frac{x}{r^{i_r}})^2 + \ldots + (\frac{x}{r^{i_r}})^{e-1}).
$$

Therefore
\begin{eqnarray*}
L_r(\bi)^{N'} & = & f(X_r^{p^{s}}) = e^{-1}(1 + \frac{X_r^{p^s}}{r^{i_r}} + (\frac{X_r^{p^s}}{r^{i_r}})^2 + \ldots + (\frac{X_r^{p^s}}{r^{i_r}})^{e-1})\\
& = & e^{-1}(1 + \frac{X_r^{p^s}}{q^{p^s{\cdot}i_r}} + (\frac{X_r^{p^s}}{q^{p^s{\cdot}i_r}})^2 + \ldots + (\frac{X_r^{p^s}}{q^{p^s{\cdot}i_r}})^{e-1})
\end{eqnarray*}
which completes the proof. \endproof

Finally we can get an explicit expression of $e(\bi)$ for the non-degenerate $H_n^\Lambda$ using \autoref{nondeg:3} and the orthogonality of $e(\bi)$'s.

\begin{Theorem} \label{nondeg:exp1}

Suppose $\bi = (i_1,i_2,\ldots,i_n) \in I^n$ and $q \neq 1$, we have
$$
e(\bi) = e^{-n} \prod_{r = 1}^n (1 + \frac{X_r^{p^s}}{q^{p^s{\cdot}i_r}} + (\frac{X_r^{p^s}}{q^{p^s{\cdot}i_r}})^2 + \ldots + (\frac{X_r^{p^s}}{q^{p^s{\cdot}i_r}})^{e-1})
$$
for $s \gg 0$.

\end{Theorem}

\section{Periodic property of $x_r$ in degenerate case}

In the degenerate cyclotomic Hecke algebras, suppose $e = p > 0$ and $H_n^\Lambda(\Z_p)$ is the degenerate cyclotomic Hecke algebra over $\Z_p$. It is obvious that $H_n^\Lambda = H_n^\Lambda(\Z_p)\otimes \F_p$. When $p > 0$ the algebra $H_n^\Lambda(\Z_p)$ is finite. We know that $\text{dim }H_n^\Lambda = \l^n n!$. Hence over $\Z_p$ the algebra has $p \l^n n!$  elements. Therefore, by choosing $k > p \l^n n!$, for any $r$ we must be able to find $k_1,k_2$ with $1 \leq k_1 < k_2 \leq k$ such that $x_r^{k_1} = x_r^{k_2}$. Therefore for any $r$ we can find integers $d_r$ and $N$ such that for any $N' \geq N$, $x_r^{N'} = x_r^{N' + d_r}$. Because $H_n^\Lambda = H_n^\Lambda(\Z_p)\otimes \F_p$, we have the same property for $x_r$'s in $H_n^\Lambda$. We define the \textbf{period} of $x_r$ to be the smallest positive integer $d_r$ such that $x_r^N = x_r^{N + d_r}$ for some $N$.

In this section we will give information on the $d_r$ and the minimal $N$ such that $x_r^N = x_r^{N+d_r}$ when $e = p > 0$.

\subsection{Form of the period $d_r$}

Recall $y_r := \sum_{\bi\in I^n} (x_r - i_r) e(\bi)$ in degenerate cyclotomic Hecke algebras by (\ref{E:y_r}).

\begin{Lemma} \label{pe:d:1}
Suppose $s$ is an integer. For any $r$, $x_r^{p^{s+1}} = x_r^{p^s}$ if and only if $y_r^{p^s} = 0$.
\end{Lemma}

\proof For any $i \in I$, we have
$$
x_r^{p^{s+1}} - x_r^{p^s} = (x_r^{p^{s+1}} - i) - (x_r^{p^s} - i) = (x_r - i)^{p^{s+1}} - (x_r - i)^{p^s}.
$$

Suppose $y_r^{p^s} = \sum_{\bi\in I^n} (x_r - i_r)^{p^s} e(\bi) = 0$. Therefore for any $\bi$, $(x_r - i_r)^{p^s} e(\bi) = 0$. Then for any $\bi \in I^n$ with $i_r = i$ we have
$$
(x_r^{p^{s+1}} - x_r^{p^s})e(\bi) = (x_r - i)^{p^{s+1}}e(\bi) - (x_r - i)^{p^s}e(\bi) = 0.
$$

Then as $\sum_{\bi\in I^n} e(\bi) = 1$, we have
$$
(x_r^{p^{s+1}} - x_r^{p^s}) = \sum_{\bi\in I^n}(x_r^{p^{s+1}} - x_r^{p^s})e(\bi) = 0,
$$
which shows that $x_r^{p^{s+1}} = x_r^{p^s}$. 

Suppose $y_r^{p^s} = \sum_{\bi\in I^n} (x_r - i_r)^{p^s} e(\bi) \neq 0$, we must be able to find a $\bi\in I^n$ such that $(x_r - i_r)^{p^s} e(\bi) \neq 0$. Assume
$$
(x_r^{p^{s+1}} - x_r^{p^s})e(\bi) = (x_r - i)^{p^{s+1}}e(\bi) - (x_r - i)^{p^s}e(\bi) = 0,
$$
which means that $y_r^{p^{s+1}} e(\bi) = y_r^{p^s}e(\bi)\neq 0$. Because $p^{s+1} > p^s$ and $y_r^{p^s}e(\bi) \neq 0$, we can find $k$ where $y_r^{p^{s+1} + k}e(\bi) = 0$ and $y_r^{p^s + k}e(\bi) \neq 0$. But $y_r^{p^{s+1} + k}e(\bi) = y_r^k y_r^{p^{s+1}}e(\bi) = y_r^k y_r^{p^s}e(\bi) = y_r^{p^s + k} e(\bi) \neq 0$, which leads to contradiction. Therefore we must have $(x_r^{p^{s+1}} - x_r^{p^s})e(\bi) = (x_r - i)^{p^{s+1}}e(\bi) - (x_r - i)^{p^s}e(\bi) \neq 0$, which yields that $x_r^{p^{s+1}} \neq x_r^{p^s}$. \endproof

Choose $s \gg 0$ such that $y_r^s = 0$. By \autoref{pe:d:1} we have $x_r^{p^s} = x_r^{p^{s+1}} = x_r^{p^s + (p-1)p^s}$. So the period $d_r$ divides $p^s(p-1)$. Then $d_r = p^m$ or $p^m(p-1)$ with $m \geq 0$.

\begin{Lemma} \label{pe:d:2}
Suppose $d_r$ is the period of $x_r$. Then $(p-1)\ |\ d_r$.
\end{Lemma}

\proof When $p = 2$ there is nothing to prove. Hence we set $p > 2$ so that $p$ is odd. Assume that $d_r = p^m$ for some $m$. Consider $\lambda = (r-1,1^{n-r+1})$ and $\t = \t^\lambda$. Let $\bj = (j_1,j_2,\ldots,j_n) = \res(\t)$, it is easy to see that $j_r = e - 1 = p-1$. Now $\bj$ is a residue sequence so that $e(\bj) \neq 0$ by \cite[Lemma 4.1(c)]{HuMathas:GradedCellular}. So we must have $\sum_{\substack{\bi\in I^n \\ i_r  = p-1}}e(\bi) \neq 0$. Choose $s\gg m$. By \autoref{deg:4},
\begin{eqnarray*}
L_r(\bj)^{N_r(\bj)} & = & -\frac{x_r^{p^s}}{p-1} - \frac{x_r^{2p^s}}{(p-1)^2} - \ldots - \frac{x_r^{(p-1)p^s}}{(p-1)^{p-1}}\\
& = & x_r^{p^s} - x_r^{2p^s} + x_r^{3p^s} - \ldots - x_r^{(p-1)p^s}.
\end{eqnarray*}

By assumption, because $s \gg m$, we have $x_r^{p^s} = x_r^{2p^s} = \ldots  = x_r^{(p-1)p^s}$. Therefore
$$
L_r(\bj)^{N_r(\bj)} = x_r^{p^s} - x_r^{2p^s} + x_r^{3p^s} - \ldots - x_r^{(p-1)p^s} = (1 - 1 + 1 - \ldots -1)x_r^{p^s} = 0.
$$

But by \autoref{deg:2} we have $L_r(\bj)^{N_r(\bj)} = \sum_{i_r = p-1} e(\bi) \neq 0$, which leads to contradiction. Therefore $d_r = p^m(p-1)$ and hence $(p-1)|d_r$. \endproof

\subsection{Period of $x_r$}

Now we know that $d_r = p^m (p-1)$ for some $m$. We can give a more specific value of $m$. Define $l$ to be the integer such that $y_r^{p^l} = 0$ and $y_r^{p^{l-1}} \neq 0$. First we introduce two Lemmas.

\begin{Lemma} \label{pe:d:3:1}
Suppose $f(x) \in \F_p[x]$, $h \in H_n^\Lambda$ and $e(\bi)h \neq 0$. Then $f(x_r) e(\bi)h = 0$ only if $f(i_r) = 0$. 
\end{Lemma}

\proof We prove this Lemma by contradiction. Because $\F_p$ is a field, $f(i_r) = 0$ only if $(x - i_r)\ |\ f(x)$. Assume $f(x_r) e(\bi)h = 0$. Suppose $f(i_r) \neq 0$, we can write $f(x) = (x - i_r)g(x) + j$ with $j \neq 0$. Set $s \gg 0$ such that $(x_r - i_r)^{p^s} e(\bi) = 0$. Because $f(x_r)e(\bi)h = 0$, we have $$
f^{p^s}(x_r)e(\bi)h = ((x_r - i_r)g(x_r) + j)^{p^s} e(\bi)h = g^{p^s}(x_r) (x_r - i_r)^{p^s} e(\bi)h + j{\cdot} e(\bi)h = j{\cdot}e(\bi)h \neq 0
$$
because $j\neq 0$ and $e(\bi)h \neq 0$, which leads to contradiction. Therefore $f(x_r)e(\bi)h \neq 0$ when $f(i_r) \neq 0$. This completes the proof.\endproof

\begin{Lemma} \label{pe:d:3:2}
Suppose $k \in \Z$ and $t \in \Z$. For any $i\in I$ with $i \neq 0$, we have
$$
x^k - x^{k + p^t(p-1)} = f(x)(i - x)^{p^t}
$$
with $f(x) = \frac{x^k}{i}(1 + \frac{x^{p^t}}{i} + (\frac{x^{p^t}}{i})^2 + \ldots + (\frac{x^{p^t}}{i})^{p-2})$.
\end{Lemma}

\proof Suppose $i \in I$ and $i \neq 0$. We have
\begin{eqnarray*}
x^k - x^{k + p^t(p-1)}
& = & x^k(1 - x^{p^t(p-1)}) = x^k(1 - (\frac{x^{p^t}}{i})^{p-1})\\
& = & x^k(1 + \frac{x^{p^t}}{i} + (\frac{x^{p^t}}{i})^2 + \ldots + (\frac{x^{p^t}}{i})^{p-2})(1 - \frac{x^{p^t}}{i})\\
& = & \frac{x^k}{i}(1 + \frac{x^{p^t}}{i} + (\frac{x^{p^t}}{i})^2 + \ldots + (\frac{x^{p^t}}{i})^{p-2})(i^{p^t} - x^{p^t})\\
& = & \frac{x^k}{i}(1 + \frac{x^{p^t}}{i} + (\frac{x^{p^t}}{i})^2 + \ldots + (\frac{x^{p^t}}{i})^{p-2})(i - x)^{p^t}\\
& = & f(x)(i - x)^{p^t}
\end{eqnarray*}
with $f(x) = \frac{x^k}{i}(1 + \frac{x^{p^t}}{i} + (\frac{x^{p^t}}{i})^2 + \ldots + (\frac{x^{p^t}}{i})^{p-2})$. This completes the proof. \endproof

\begin{Proposition} \label{pe:d:3}
Suppse $l$ is the smallest non-negative integer such that $y_r^{p^l} = 0$. Then the period of $x_r$ is $d_r = p^l(p-1)$.
\end{Proposition}

\proof Suppose $d_r = p^m (p-1)$. By \autoref{pe:d:1} we have $x_r^{p^{l+1}} = x_r^{p^l + p^l(p-1)} = x_r^{p^l}$. Therefore $d_r\ |\ p^l(p-1)$ which indicates that $m \leq l$. Now take $s \gg 0$, by \autoref{pe:d:3:2} we have
$$
(x_r^{p^s} - x_r^{p^s + p^{l-1}(p-1)})e(\bi) = f(x_r) e(\bi) (i_r - x_r)^{p^{l-1}} 
$$
where $f(x) = \frac{x^{p^s}}{i_r}(1 + \frac{x^{p^{l-1}}}{i_r} + (\frac{x^{p^{l-1}}}{i_r})^2 + \ldots + (\frac{x^{p^{l-1}}}{i_r})^{p-2}) \in \F_p[x]$. It is easy to see that $f(i_r) = p-1 \neq 0$. By the definition of $l$, $e(\bi)(i_r - x_r)^{p^{l-1}} \neq 0$. Then by \autoref{pe:d:3:1} we have
$$
(x_r^{p^s} - x_r^{p^s + p^{l-1}(p-1)})e(\bi) = f(x_r) e(\bi) (i_r - x_r)^{p^{l-1}} \neq 0.
$$

Therefore $x_r^{p^s} - x_r^{p^s + p^{l-1}(p-1)} \neq 0$, i.e. $x_r^{p^s} \neq x_r^{p^s + p^{l-1}(p-1)}$, which yields $m \geq l$. This shows that $ m = l$ and $d_r = p^l(p-1)$. \endproof

\subsection{Periodic property of $x_r$}

Now we know that the period of $x_r$ is $d_r = p^l(p-1)$, and we still need to find the smallest non-negative integer $N$ such that $x_r^N = x_r^{N+d_r}$.

\begin{Proposition} \label{pe:d:4}
Suppose $1\leq r\leq n$ and we can find a residue sequence $\bi$ such that $i_r = 0$. If $N$ is the smallest non-negative integer such that $x_r^N \sum_{i_r = 0}e(\bi) = 0$, then $x_r^N = x_r^{N+d_r}$ and $x_r^{N-1} \neq x_r^{N-1+d_r}$.
\end{Proposition}

\proof By the definition of $N$, we can find $\bi$ with $i_r = 0$ such that $x_r^{N-1}e(\bi) \neq 0$ and $x_r^N e(\bi) = 0$. Suppose $s \gg 0$. Because $d_r \geq 1$ we have
$$
(x_r^{N-1} - x_r^{N-1+d_r})e(\bi) = x_r^{N-1} e(\bi) - x_r^{N-1+d_r} e(\bi) = x_r^{N-1} e(\bi) \neq 0
$$
which indicates that $x_r^{N-1} \neq x_r^{N-1+d_r}$.

Next we will prove that $x_r^N = x_r^{N+d_r}$. Suppose $\bi\in I^n$ with $i_r = 0$, then
$$
(x_r^N - x_r^{N+d_r})e(\bi) = (1 - x_r^{d_r})x_r^N e(\bi) = 0
$$
by the definition of $N$. Now suppose $\bi \in I^n$ with $i_r \neq 0$. By \autoref{pe:d:3}, $d_r = p^l(p-1)$ where $y_r^{p^l} = 0$.  So by \autoref{pe:d:3:2},
$$
(x_r^N - x_r^{N+d_r})e(\bi) = (x_r^N - x_r^{N+(p-1)p^l})e(\bi) = f(x_r)e(\bi)(i_r - x_r)^{p^l} = f(x_r)e(\bi)(-y_r)^{p^l} = 0
$$
with $f(x) \in \F_p[x]$. Therefore we have $(x_r^N - x_r^{N+d_r})e(\bi) = 0$ for any $\bi \in I^n$ and hence $x_r^N = x_r^{N+d_r}$. This completes the proof. \endproof

Notice that in \autoref{pe:d:4} we require $1 \leq r \leq n$ such that we can find a residue sequence $\bi$ with $i_r = 0$. If no such residue sequence exists we obtain a different result.

\begin{Proposition} \label{pe:d:5}
Suppose $1\leq r\leq n$ and for any residue sequence $\bi$ we always have $i_r \neq 0$. Then $x_r^{d_r} = 1$.
\end{Proposition}

\proof By \autoref{pe:d:3}, $d_r = p^l (p-1)$ where $y_r^{p^l} = 0$. And for any $\bi \in I^n$, we have $i_r \neq 0$. Then by \autoref{pe:d:3:2},
$$
(1 - x_r^{d_r})e(\bi) = (1 - x_r^{p^l(p-1)}) e(\bi) = f(x_r) e(\bi) (i_r - x_r)^{p^l} = f(x_r)e(\bi)(-y_r)^{p^l} = 0,
$$
which shows that $x_r^{d_r}e(\bi) = e(\bi)$ for any $\bi\in I^n$. Hence $x_r^{d_r} = 1$.\endproof

Finally we give the main Theorem of this section by combining \autoref{pe:d:3}, \autoref{pe:d:4} and \autoref{pe:d:5}.

\begin{Theorem} \label{pe:d}
In the degenerate cyclotomic Hecke algebras, suppose $l$ is the smallest nonnegative integer such that $y_r^{p^l} = 0$ and $N$ is the smallest nonnegative integer such that $x_r^N\sum_{i_r = 0}e(\bi) = 0$. Then $x_r^k = x_r^{k + p^m(p-1)}$ if and only if $m \geq l$ and $k\geq N$.
\end{Theorem}

\section{Periodic property of $X_r$ in non-degenerate case}

In the non-degenerate cyclotomic Hecke algebras, suppose $H_n^\Lambda(\Z_p[q])$ is the non-degenerate cyclotomic Hecke algebra over $\Z_p[q]$. When $e > 0$ and $p > 0$, $\Z_p[q]$ is finite. Therefore for the same reason as in degenerate cases, we can find $N$ and $d_r$ such that $X_r^N = X_r^{N+d_r}$, i.e. $X_r$ has periodic property in $H_n^\Lambda(\Z_p[q])$. Because $H_n^\Lambda= H_n^\Lambda(\Z_p[q])\otimes \F_p$, we have the same property for $X_r$'s in $H_n^\Lambda$. The period $d_r$ in non-degenerate cases is defined similarly as in degenerate cases.

In this section we will give information on the $d_r$ and the minimal $N$ such that $X_r^N = X_r^{N+d_r}$ when $e > 0$ and $p > 0$.

\subsection{Form of period $d_r$}

Recall that $y_r = \sum_{\bi\in I^n} (1 - q^{-i_r}X_r) e(\bi)$ in non-degenerate cyclotomic Hecke algebras by (\ref{E:y_r}).

\begin{Lemma} \label{pe:n:1}
Suppose $s \gg 0$ and $1 \leq r \leq n$. We have $X_r^{ep^s} = 1$.
\end{Lemma}

\proof By \autoref{nondeg:3}, for any $i_r \in I$, we have
\begin{eqnarray*}
&&(X_r - q^{i_r})^{p^s}\sum_{\substack{\bj \in I^n \\ j_r = i_r}} e(\bj) =  (X_r^{p^s} - q^{p^s{\cdot}i_r})\sum_{\substack{\bj \in I^n \\ j_r = i_r}} e(\bj)\\
& = & e^{-1}(X_r^{p^s} - q^{p^s{\cdot}i_r})(1 + \frac{X_r^{p^s}}{q^{p^s{\cdot}i_r}} + (\frac{X_r^{p^s}}{q^{p^s{\cdot}i_r}})^2 + \ldots + (\frac{X_r^{p^s}}{q^{p^s{\cdot}i_r}})^{e-1})\\
& = & e^{-1}(\frac{X_r^{ep^s}}{q^{(e-1)p^s{\cdot}i_r}} - q^{p^s{\cdot}i_r}) = 0,
\end{eqnarray*}
which leads to
$$
\frac{X_r^{ep^s}}{q^{(e-1)p^s{\cdot}i_r}} = q^{p^s{\cdot}i_r} \hspace*{3mm} \Rightarrow \hspace*{3mm} X_r^{ep^s} = q^{ep^s{\cdot}i_r} = 1
$$
because $q^e = 1$. \endproof

Define $d_r$ to be the period of $X_r$. By \autoref{pe:n:1} we have $d_r\ |\ ep^s$. Therefore $d_r = ep^m$ with $m \geq 0$ or $d_r = p^m$ with $m \geq 1$. In the following Lemma we are going to give more information about the form of $d_r$.

\begin{Lemma} \label{pe:n:2}
Suppose $d_1$ is the period of $X_1$. We have $d_1 = p^m$ if $\Lambda = \l \Lambda_0$.
\end{Lemma}

\proof By (\ref{E:cycloHecke}) we have $(X_1 - q^0)^\l = (X_1 - 1)^\l = 0$. Choose $s$ such that $p^s \geq \l$, we have $(X_1 - 1)^{p^s} = X_1^{p^s} - 1 = 0$, which means $X_1^{p^s} = 1$. Hence $d_1\ |\ p^s$ and therefore $d_1 = p^m$.\endproof

\begin{Remark}
When we set $r = 1$ and $\Lambda = \l \Lambda_0$, it means that $e(\bi) = 0$ if $i_r = i_1 \neq 0$. So \autoref{pe:n:2} is actually:\\

\textit{Suppose $1 \leq r \leq n$ and for any $\bi \in I^n$, $e(\bi) = 0$ if $i_r \neq 0$. Then $d_r = p^m$.}\\

\hspace*{-4mm}because the only possible $r$ and $\Lambda$ for such condition is giving in \autoref{pe:n:2}.

\end{Remark}

\begin{Lemma} \label{pe:n:3}
Suppose $d_r$ is the period of $X_r$. We have $e\ |\ d_r$ if $r > 1$ or $r = 1$ and $\Lambda \neq \l\Lambda_0$.
\end{Lemma}

\proof We prove the Lemma by contradiction. Assume that $d_r = p^m$. Choose $i_r \in I$ with $i_r \neq 0$. Because $r > 1$ or $r = 1$ and $\Lambda \neq \l\Lambda_0$ we must can find $\bj \in I^n$ with $j_r = i_r$ with $e(\bj) \neq 0$. Then $\sum_{\substack{\bj\in I^n \\j_r = i_r}}e(\bj) \neq 0$. Choose $s \gg m$. By \autoref{e,p in nondeg}, $\text{gcd}(e,p) = 1$. Then because $i_r \neq 0$, $p^s{\cdot}i_r \not\equiv 0\pmod{e}$. Then by \autoref{nondeg:1} and \autoref{nondeg:3}, we have
\begin{eqnarray*}
\sum_{\substack{\bj\in I^n \\j_r = i_r}}e(\bj) & = & e^{-1}(1 + \frac{X_r^{p^s}}{q^{p^s{\cdot}i_r}} + (\frac{X_r^{p^s}}{q^{p^s{\cdot}i_r}})^2 + \ldots + (\frac{X_r^{p^s}}{q^{p^s{\cdot}i_r}})^{e-1})\\
& = & e^{-1}(1 + \frac{1}{q^{p^s{\cdot}i_r}} + \frac{1}{q^{2p^s{\cdot}i_r}} + \ldots + \frac{1}{q^{(e-1)p^s{\cdot}i_r}})X_r^{p^s}\\
& = & \frac{e^{-1}}{q^{(e-1)p^s{\cdot}i_r}}(1 + q^{p^s{\cdot}i_r} + (q^{p^s{\cdot}i_r})^2 + \ldots + (q^{p^s{\cdot}i_r})^{e-1})X_r^{p^s} = 0,
\end{eqnarray*}
which leads to contradiction. Hence $d_r \neq p^m$ and therefore $e\ |\ d_r$. \endproof

\subsection{Period of $X_r$}
Now we know that $d_r = p^m$ when $r = 1$ and $\Lambda = \l\Lambda_0$ and $d_r = ep^m$ otherwise. In the rest of the section we will find the value of $m$. First we give the simpler case.

\begin{Lemma} \label{pe:n:4}
Suppose $s \geq 0$ and $\Lambda = \l \Lambda_0$. We have $X_1^{p^s} = 1$ if and only if $y_1^{p^s} = 0$.
\end{Lemma}

\proof Suppose $y_1^{p^s} = 0$. For any $\bi \in I^n$,
$$
y_1^{p^s}e(\bi) = (1 - X_1)^{p^s}e(\bi) = e(\bi) - X_1^{p^s}e(\bi) = 0 \hspace*{3mm} \Rightarrow \hspace*{3mm} X_1^{p^s}e(\bi) = e(\bi).
$$

Therefore $X_1^{p^s} = 1$.

Suppose $y_1^{p^s} \neq 0$. Then we can find $\bi\in I^n$ with $y_1^{p^s}e(\bi) \neq 0$. So
$$
y_1^{p^s}e(\bi) = (1 - X_1)^{p^s}e(\bi) = e(\bi) - X_1^{p^s}e(\bi) \neq 0 \hspace*{3mm} \Rightarrow \hspace*{3mm} X_1^{p^s}e(\bi) \neq e(\bi).
$$

Therefore $X_1^{p^s} \neq 1$. \endproof

Now we consider the case when $r \neq 1$ or $r = 1$ and $\Lambda \neq \l\Lambda_0$.

\begin{Lemma} \label{q^pk = q^ps}
For any non-negative integer $s$, we can find $k \gg s$ such that $q^{p^k} = q^{p^s}$ and $p^{k - s} \equiv 1\pmod{e}$ .
\end{Lemma}

\proof By \autoref{e,p in nondeg}, we can find $l$ such that $p^l \equiv 1\pmod{e}$. Because $q^e = 1$, choose $t \gg 0$ and set $k = s + tl$, we have $q^{p^k} = q^{p^{s + tl}} = q^{p^s p^{tl}} = (q^{p^{tl}})^{p^s} = q^{p^s}$. Moreover, $p^{k-s} = p^{tl} \equiv 1^t\pmod{e} \equiv 1\pmod{e}$. This completes the proof.  \endproof

Now we are ready to give more information of $d_r$.

\begin{Lemma}\label{pe:n:5:1}
We have $X_r^{ep^s} = 1$ for some $s$ only if $y_r^{p^s} =  0$.
\end{Lemma}

\proof Fix $s$ such that $X_r^{ep^s} = 1$. By \autoref{q^pk = q^ps}, we can find $k\gg 0$ such that $q^{p^{s+k}} = q^{p^s}$ and $p^k \equiv 1 \pmod{e}$. Therefore $p^{s+k} - p^s = p^s(p^k - 1)$ and hence $ep^s\ |\ p^s(p^k-1)$. So $X_r^{p^{s+k}} = X_r^{p^s}$.

Then for any $i_r \in I$, by \autoref{nondeg:3}, we have
\begin{eqnarray*}
\sum_{\substack{\bj \in I^n \\ j_r = i_r}}e(\bj) & = & e^{-1}(1 + \frac{X_r^{p^{s+k}}}{q^{p^{s+k}{\cdot}i_r}} + (\frac{X_r^{p^{s+k}}}{q^{p^{s+k}{\cdot}i_r}})^2 + \ldots + (\frac{X_r^{p^{s+k}}}{q^{p^{s+k}{\cdot}i_r}})^{e-1})\\
& = & e^{-1}(1 + \frac{X_r^{p^s}}{q^{p^s{\cdot}i_r}} + (\frac{X_r^{p^s}}{q^{p^s{\cdot}i_r}})^2 + \ldots + (\frac{X_r^{p^s}}{q^{p^s{\cdot}i_r}})^{e-1}).
\end{eqnarray*}

Therefore,
\begin{eqnarray*}
(X_r - q^{i_r})^{p^s}\sum_{\substack{\bj \in I^n \\ j_r = i_r}}e(\bj) & = & e^{-1}(X_r - q^{i_r})^{p^s}(1 + \frac{X_r^{p^s}}{q^{p^s{\cdot}i_r}} + (\frac{X_r^{p^s}}{q^{p^s{\cdot}i_r}})^2 + \ldots + (\frac{X_r^{p^s}}{q^{p^s{\cdot}i_r}})^{e-1})\\
& = & e^{-1}(X_r^{p^s} - q^{p^s{\cdot}i_r})(1 + \frac{X_r^{p^s}}{q^{p^s{\cdot}i_r}} + (\frac{X_r^{p^s}}{q^{p^s{\cdot}i_r}})^2 + \ldots + (\frac{X_r^{p^s}}{q^{p^s{\cdot}i_r}})^{e-1})\\
& = & e^{-1} (\frac{X_r^{ep^s}}{q^{(e-1)p^s{\cdot}i_r}} - q^{p^s{\cdot}i_r}) = e^{-1} (\frac{1}{q^{(e-1)p^s{\cdot}i_r}} - q^{p^s{\cdot}i_r}) = 0
\end{eqnarray*}
because $\frac{1}{q^{(e-1)p^s{\cdot}i_r}} = q^{p^s{\cdot}i_r}$. This means that $y_r^{p^s} \sum_{\substack{\bj \in I^n \\ j_r = i_r}}e(\bj) = (1 - q^{-i_r}X_r)^{p^s}\sum_{\substack{\bj \in I^n \\ j_r = i_r}}e(\bj) = 0$ for any $i_r \in I$. Hence $y_r^{p^s} = 0$. \endproof

\begin{Lemma} \label{pe:n:5:2}
Suppose $y_r^{p^s} = 0$ for some $s$. Then we have $X_r^{ep^s} = 1$.
\end{Lemma}

\proof Fix $s$ such that $y_r^{p^s} = 0$. For any $\bi \in I^n$, by \autoref{q^pk = q^ps} we can choose $k \gg s$ such that $q^{p^k} = q^{p^s}$. Then
\begin{eqnarray*}
(X_r^{p^k} - X_r^{p^s})e(\bi) & = & (X_r^{p^k} - q^{p^s{\cdot}i_r} - X_r^{p^s} + q^{p^s{\cdot}i_r})e(\bi)\\
& = & (X_r^{p^k} - q^{p^k{\cdot}i_r})e(\bi) - (X_r^{p^s} - q^{p^s{\cdot}i_r})e(\bi)\\
& = & (X_r - q^{i_r})^{p^k}e(\bi) - (X_r - q^{i_r})^{p^s}e(\bi) = q^{p^k{\cdot}i_r}(-y_r)^{p^k}e(\bi) - q^{p^s{\cdot}i_r}(-y_r)^{p^s}e(\bi) = 0.
\end{eqnarray*}

So $(X_r^{p^k} - X_r^{p^s})e(\bi) = 0$ for any $\bi \in I^m$. Therefore we must have $X_r^{p^k} - X_r^{p^s} = 0$ for some $k\gg 0$. Hence
$$
X_r^{p^k} - X_r^{p^s} = X_r^{p^s}(X_r^{p^k - p^s} - 1) = 0 \hspace*{5mm} \Rightarrow \hspace*{5mm} X_r^{p^k - p^s} = 1,
$$
which implies that $d_r\ |\ p^k - p^s$. We know that $d_r = ep^m$ for some $m$ and $p^k - p^s = p^s(p^{k-s}-1)$. It is obvious that $m \leq s$. Hence $X_r^{ep^s} = 1$. \endproof

The next Corollary follows straightforward by combining \autoref{pe:n:5:1} and \autoref{pe:n:5:2}.

\begin{Corollary} \label{pe:n:5}
Suppose $s\geq 0$, $r > 1$ or $r = 1$ and $\Lambda \neq \l\Lambda_0$. We have $X_r^{ep^s} = 1$ if and only if $y_r^{p^s} = 0$.
\end{Corollary}

Finally, combining all the results above, we have the final Theorem.

\begin{Theorem} \label{pe:n}
In non-degenerate $H_n^\Lambda$, we have $X_r^{d_r} = 1$ with
$$
d_r = \begin{cases}
p^m, & \text{if $r = 1$ and $\Lambda = \l\Lambda_0$;}\\
ep^m, & \text{otherwise.}
\end{cases}
$$
if and only if $y_r^{p^m} = 0$.
\end{Theorem}

\proof The Theorem follows straightforward by \autoref{pe:n:2}, \autoref{pe:n:3}, \autoref{pe:n:4} and \autoref{pe:n:5}. \endproof

\section{Explicit expression of idempotents $e(\bi)$ part II}

In \autoref{deg:exp1} and \autoref{nondeg:exp1}, we gave explicit expressions of $e(\bi)$'s in degenerate and non-degenerate cyclotomic Hecke algebras. But the expression depends on the choice of $s$, and up to now the only information we know about $s$ is that $s \gg 0$. In this section we give a more specific value for $s$ in both degenerate and non-degenerate case.

Define $l_r$ to be the smallest integer such that $y_r^{p^{l_r}} = 0$.

\begin{Theorem}\label{deg:exp2}
Suppose $\bi = (i_1,i_2,\ldots,i_n) \in I^n$ and $q = 1$. Then $e(\bi) = \prod_{r = 1}^n \Lr$, where
$$
\Lr  = \begin{cases}
1 - x_r^{p^{l_r}(1-p)}, & \text{when $i_r = 0$};\\
-\sum_{k = 1}^{p-1} \frac{x_r^{kp^{l_r}}}{i_r^k}, & \text{when $i_r \neq 0$}.
\end{cases}
$$
\end{Theorem}

\proof By \autoref{deg:exp1}, we have that $e(\bi) = \prod_{r = 1}^n \Lr$, where
$$
\Lr  = \begin{cases}
1 - x_r^{p^s(1-p)}, & \text{when $i_r = 0$};\\
-\sum_{k = 1}^{p-1} \frac{x_r^{kp^s}}{i_r^k}, & \text{when $i_r \neq 0$}.
\end{cases}
$$
for $s \gg 0$. Then by \autoref{pe:d:1}, it is easy to see that when $s \gg 0$, $x_r^{p^s} = x_r^{p^{l_r}}$, and hence
\begin{align*}
1 - x_r^{p^s(1-p)} &= 1 - x_r^{p^{l_r}(1-p)};\\
-\sum_{k = 1}^{p-1} \frac{x_r^{kp^s}}{i_r^k} & = -\sum_{k = 1}^{p-1} \frac{x_r^{kp^{l_r}}}{i_r^k},
\end{align*}
which completes the proof. \endproof

\begin{Theorem} \label{nondeg:exp2}

Suppose $\bi = (i_1,i_2,\ldots,i_n) \in I^n$ and $q \neq 1$, we have
$$
e(\bi) = e^{-n} \prod_{r = 1}^n (1 + \frac{X_r^{p^{l_r}}}{q^{p^{l_r}{\cdot}i_r}} + (\frac{X_r^{p^{l_r}}}{q^{p^{l_r}{\cdot}i_r}})^2 + \ldots + (\frac{X_r^{p^{l_r}}}{q^{p^{l_r}{\cdot}i_r}})^{e-1}).
$$

\end{Theorem}

\proof Suppose $\bi = (i_1,i_2,\ldots,i_n) \in I^n$ and $q \neq 1$. By \autoref{nondeg:exp1}, we have
$$
e(\bi) = e^{-n} \prod_{r = 1}^n (1 + \frac{X_r^{p^s}}{q^{p^s{\cdot}i_r}} + (\frac{X_r^{p^s}}{q^{p^s{\cdot}i_r}})^2 + \ldots + (\frac{X_r^{p^s}}{q^{p^s{\cdot}i_r}})^{e-1})
$$
for $s \gg 0$. If $\Lambda = \l \Lambda_0$, for any $\bi \in I^n$ with $e(\bi) \neq 0$, by \cite[Lemma 4.1(c)]{HuMathas:GradedCellular} we must have $i_1 = 0$.

Then we have
$$
e(\bi) = e^{-n} \prod_{r = 2}^n (1 + \frac{X_r^{p^s}}{q^{p^s{\cdot}i_r}} + (\frac{X_r^{p^s}}{q^{p^s{\cdot}i_r}})^2 + \ldots + (\frac{X_r^{p^s}}{q^{p^s{\cdot}i_r}})^{e-1})
$$
with $s \gg 0$. By \autoref{q^pk = q^ps}, choose $s \gg l_r$ such that $q^{p^s} = q^{p^{l_r}}$ and $p^{s - l_r} \equiv 1\pmod{e}$. Because $p^s - p^{l_r} = (p^{s-l_r}-1)p^{l_r}$, $ep^{l_r}\ |\ (p^s - p^l)$. Therefore by \autoref{pe:n}, for any $r \geq 2$, $X_r^{p^s} = X_r^{p^{l_r} + p^s - p^{l_r}} = X_r^{p^{l_r}}$. This means that
$$
1 + \frac{X_r^{p^s}}{q^{p^s{\cdot}i_r}} + (\frac{X_r^{p^s}}{q^{p^s{\cdot}i_r}})^2 + \ldots + (\frac{X_r^{p^s}}{q^{p^s{\cdot}i_r}})^{e-1} = 1 + \frac{X_r^{p^{l_r}}}{q^{p^{l_r}{\cdot}i_r}} + (\frac{X_r^{p^{l_r}}}{q^{p^{l_r}{\cdot}i_r}})^2 + \ldots + (\frac{X_r^{p^{l_r}}}{q^{p^{l_r}{\cdot}i_r}})^{e-1}
$$
for any $r \geq 2$. Hence
$$
e(\bi) = e^{-n} \prod_{r = 2}^n (1 + \frac{X_r^{p^{l_r}}}{q^{p^{l_r}{\cdot}i_r}} + (\frac{X_r^{p^{l_r}}}{q^{p^{l_r}{\cdot}i_r}})^2 + \ldots + (\frac{X_r^{p^{l_r}}}{q^{p^{l_r}{\cdot}i_r}})^{e-1}).
$$

Now recall $y_r^{p^l} = 0$, we have $X_r^{p^l}e(\bi) = 0$ because $i_1 = 0$. Therefore
$$
e(\bi)(1 + \frac{X_r^{p^{l_r}}}{q^{p^{l_r}{\cdot}i_r}} + (\frac{X_r^{p^{l_r}}}{q^{p^{l_r}{\cdot}i_r}})^2 + \ldots + (\frac{X_r^{p^{l_r}}}{q^{p^{l_r}{\cdot}i_r}})^{e-1}) = e(\bi),
$$
and hence
$$
e(\bi) = e^{-n} \prod_{r = 1}^n (1 + \frac{X_r^{p^{l_r}}}{q^{p^{l_r}{\cdot}i_r}} + (\frac{X_r^{p^{l_r}}}{q^{p^{l_r}{\cdot}i_r}})^2 + \ldots + (\frac{X_r^{p^{l_r}}}{q^{p^{l_r}{\cdot}i_r}})^{e-1}).
$$

If $\Lambda \neq \l\Lambda_0$, by \autoref{q^pk = q^ps}, choose $s \gg l$ such that $q^{p^s} = q^{p^l}$ and $p^{s - l} \equiv 1\pmod{e}$. Hence for the same reason as when $\Lambda = \l\Lambda_0$, we have
$$
1 + \frac{X_r^{p^s}}{q^{p^s{\cdot}i_r}} + (\frac{X_r^{p^s}}{q^{p^s{\cdot}i_r}})^2 + \ldots + (\frac{X_r^{p^s}}{q^{p^s{\cdot}i_r}})^{e-1} = 1 + \frac{X_r^{p^{l_r}}}{q^{p^{l_r}{\cdot}i_r}} + (\frac{X_r^{p^{l_r}}}{q^{p^{l_r}{\cdot}i_r}})^2 + \ldots + (\frac{X_r^{p^{l_r}}}{q^{p^{l_r}{\cdot}i_r}})^{e-1}
$$
for any $r$. Therefore
$$
e(\bi) = e^{-n} \prod_{r = 1}^n (1 + \frac{X_r^{p^{l_r}}}{q^{p^{l_r}{\cdot}i_r}} + (\frac{X_r^{p^{l_r}}}{q^{p^{l_r}{\cdot}i_r}})^2 + \ldots + (\frac{X_r^{p^{l_r}}}{q^{p^{l_r}{\cdot}i_r}})^{e-1}).
$$
\endproof

We remark that we do not know explicit expressions for the $e(\bi)$'s and in general the $x_r$ and $X_r$'s do not have a periodic property when $p = 0$ or $e = 0$. The main problem is because we do not have analogue result of~\autoref{useful} when $p = 0$ and~\autoref{nondeg:2} also fails when $e = 0$. But~\autoref{deg:1} is still true when $p = 0$ or $e = 0$. Define $I_r(\bi) := \set{j \in I|j \neq i_r, \text{ there exists $e(\bj) \neq 0$ with $j_r = j$}}$ and set $L'_r(\bi) = \prod_{j \in I_r(\bi)} L_{i_r,j}$. By \cite[Lemma 4.1(c)]{HuMathas:GradedCellular} we have $\set{e(\bi)|e(\bi) \neq 0}$ is finite and hence $L'_r(\bi)$ is well-defined. Then $L'_r(\bi)$ satisfy \autoref{deg:2} and \autoref{deg:3}. Hence even though we cannot give a simplified expression of $e(\bi)$'s similar to~\autoref{deg:exp2} and~\autoref{nondeg:exp2} when $p = 0$ or $e = 0$, we can still have an expression for $e(\bi)$'s which is useful for computation purposes.

\bibliography{papers}

\end{document}